# Two-stage robust bilevel optimization model for facility location considering operational service level under disruption risk

Haitao HU

School of Management Science & Engineering, Dongbei University of Finance and Economics

haitaohudufe@163.com

Jing LI

School of Management Science & Engineering, Dongbei University of Finance and Economics

lijing2020dufe@163.com

Jiafu TANG

School of Management Science & Engineering, Dongbei University of Finance and Economics

tangjiafu@dufe.edu.cn

Bo ZENG

Department of Industrial Engineering, University of Pittsburgh

bzeng@pitt.edu

**Abstract**

The bilevel facility location problem (BO-FLP) is one of the core optimization problems behind the design of many decentralized industrial systems, e.g., supply chain systems where a supplier constructs some critical facilities and then uses them to serve retailers in a cost-effective fashion, while retailers directly handle customers aiming to minimize the total unmet demand in a rather independent fashion. When uncertainty is considered, scenario-based stochastic approaches are commonly used, but they often become impractical due to insufficient data or an exponential number of scenarios. To address this issue, this paper adopts robust optimization and proposes a novel two-stage robust bilevel facility location model. Several structural properties are derived to improve both theoretical understanding and solution efficiency. Based on this, an enhanced column-and-constraint generation algorithm is developed for robust bilevel optimization with decision-dependent uncertainty, significantly improving exact solution capability over the standard method. Numerical results show that, compared to the centralized two-stage RO model, our model pays more attention to demand fulfillment, typically resulting in higher service efficiency and better utilization of supply capacity. Under a small-scale disruption, this new model delivers better service performance. However, under a large-scale disruption, the centralized model performs more effectively.

## 1. Introduction

The facility location problem (FLP) is a fundamental topic in operations research, which is often used to address the design and operational problems in supply chain, logistics and telecommunications systems. In simple and most common situations, it is often formulated as a mixed integer program (MIP), which seeks to minimize the overall cost that includes both the fixed cost associated with facility construction decisions and the operational cost associated with service decisions between facilities and clients (Snyder et al., 2016). We note that this single-level formulation depends on a key assumption — namely, that all decisions are made by a single decision maker (DM) acting in pursuit of a unified interest. In other words, the model presumes a vertical decision-making system where no conflict exists among decision-making parties, as shown in Fig. 1a.

Nevertheless, very often this is not the case in supply chain systems. Many real-world systems, decisions are made by multiple stakeholders who may have distinct, and sometimes conflicting, interests. For instance, a DM (e.g., a manufacturer), may determine the locations of facilities to optimize her overall cost or efficiency, acting as the upper-level DM that is also known as the leader. While individual retailers or regional managers in the system make local decisions, e.g., order quantities or service levels, based on their own interests or specifications, acting as the lower-level DMs that is also known as the followers. Under such circumstances, the whole supply chain system is decentralized and hierarchical, as illustrated in Fig. 1b, which is generally modeled by the bilevel optimization (BO) scheme.

For instance, in the context of distribution planning, Camacho-Vallejo et al. (2022) proposed a BO model in which a central firm seeks to maximize its profit while considering the transportation cost minimization behavior of downstream distributors. Similarly, in the domain of emergency facility location and disaster management, Liu et al. (2023) considered a BO model where a government agency, acting as the leader, aims to improve the performance of the trauma care network while limiting public subsidies for upgrading emergency centers, whereas hospital groups, as the follower, seek to redesign the network to maximize marginal profit. The layout of charging stations or distributed energy sites also faces similar problems. Even for ordinary electric vehicle charging piles facility location, while the system operator determining the layout, the issue of the consumers' preferences also needs to be considered by a BO model (Lamontagne et al., 2023). It is often recognized that BO provides a structured framework that enables more adaptive and resilient supply chain network design.

When random factors, e.g., disruptions on facilities, are considered, strategic and operational decisions are generally considered in two different stages, i.e., strategic ones in stage 1 and operational ones in stage 2, leading to a two-stage decision making problem. The classical MIP models have been extended to their two-stage stochastic programming (SP) or robust optimization (RO), according to the available data and attitude towards randomness (Snyder et al., 2016). For BO models, a few SP extensions have been developed and investigated (Zhang & Özaltın, 2021; Li et al, 2023). Actually, different from the centralized MIP models, the lower-level problem often coincides with the second-stage decision making, rendering the leader being the first-stage and the follower the second-stage decision makers respectively. As in Fig. 1b, the network designer, i.e., the leader, makes pre-disruption strategic location decisions, and the network user, serving as the follower, optimizes allocation and ensures supply security after the disruption risk occurs. We mention that, due to limited historical disruption data in practice, it is often difficult to accurately obtain or estimate the probability of the facility disruptions. In such cases, RO is more appropriate as it does not depend on probabilistic distribution information. Now, a few two-stage RO for centralized MIP have appeared in the literature uncertainties (An et al., 2014, Cheng et al., 2021; Ahmadi et al., 2023). However, they overlook the distinct objective of the operational level, leading to allocation decisions that may fail to meet the actual needs of the network user, and consequently, result in solutions that are misaligned with practical realities.

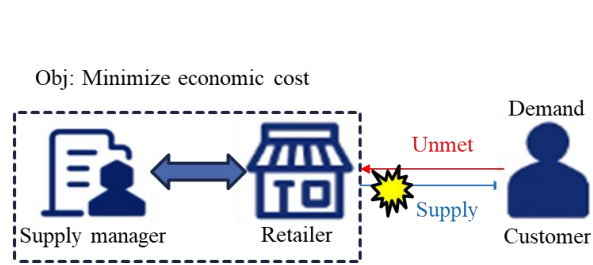
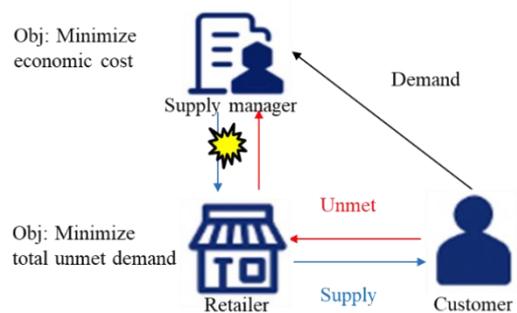

Fig.1a Centralized decision-making          Fig.1b Decentralized decision-making

Fig.1   The difference between the Centralized and Decentralized decision-making

To simultaneously address both the disruption risk and the relevant shareholders decision in real life, this paper proposes a more comprehensive model, named as two-stage robust bilevel facility location problem (2-stage RBO-FLP), that integrates the distinct features of both the RO-FLP and BO-FLP models. The 2-stage RBO-FLP model captures both the disruption risk and the behavioral divergence between the designer and the user in the FLP, providing a more realistic and adaptive approach to resilient supply chain network. Specifically, in 2-stage RBO-FLP, the upper-level network designer makes pre-disruption facility location decisions to minimize the total system cost in the first-stage. At the lower level (also the second stage) the network user responds by making optimal supply allocation decisions based on the available operational facilities after the disruption risk, aiming at minimize total unmet demand to guarantee service level. Unlike traditional two-stage RO, it features a BO formulation in the second stage, reflecting the decentralized interaction between the network designer and the user. Compared existing models, our 2-stage RBO-FLP model is significantly more challenging to solve due to its complex structure.

To solve the 2-stage RBO-FPL model exactly, we proposed an improved version of C&CG-BO, a method for general two-stage robust BO problems (Xu, 2023), by adopting decision-dependent uncertainty (DDU) to reduce the search space of disruption scenario for the 2-stage RBO-FLP model, named as C&CG-BO-DDU. Similar to the basic C&CG and C&CG-BO, C&CG-BO-DDU involves solving a master problem (MP) and a subproblem (SP) in an iterative fashion. It is worth highlighting that SP of C&CG-BO is actually a pessimistic bilevel optimization model that is very challenging to solve. Nevertheless, by leveraging reformulation techniques and the critical DDU structure of the 2-stage RBO-FLP, C&CG-BO-DDU is able to improve our solution capacity to SP drastically, and therefore leading to a significant reduction in the overall computational time.

The main contributions of this paper are summarized as follows: 1) We extend the bilevel facility location problem into uncertain scenarios under two-stage RO framework and develop a two-stage robust bilevel facility location model considering operational service level under the facility disruption risk. To the best of our knowledge, it is the first work to propose a two-stage robust BO model for facility location problem under disruption risk, which is significant for managers to pay more attention on operational performance when designing a resilient supply chain in practice. 2) We then derive several theoretical properties to highlight the fundamental differences between two-stage RO-FLP and two-stage RBO-FLP, and induce the condition under which they are same. C&CG-BO algorithm is customized and improved, referred to as C&CG-BO-DDU, that solves the 2-stage RBO-FLP model more efficiently compared to its standard implementation. 3) We conduct extensive numerical tests to highlight that explicitly integrating the operational performance into strategic supply chain network design not only reduces available supply capacity waste and improves overall demand satisfaction, but also uncovers a dynamic coordination mechanism. It shows that, under lower disruption risk, a decentralized decision-making structure can improve service level by allowing operational flexibility. In contrast, under higher disruption risk, a centralized management approach is more effective in enhancing the overall resilience and performance of the supply chain network.

The rest of this paper is organized as follows. Section 2 reviews related literature. Section 3 presents the model formulation and deduces some properties. Section 4 develops an improved C&CG-BO and presents our enhancement. Section 5 presents the numerical experiments and results analysis. Section 6 concludes the paper.

**2. Literature review**

Focusing on the facility location studied in this paper, there are two streams in the literature. One is single level optimization and BO from perspectives of decision process, i.e., centralization or decentralization, and the other is SP or RO models, from the perspectives of modeling of uncertainty. In the follows, we will review the state-of-the-art of the above models.

2.1 Single level approaches for facility location

For single level FLP, its SP and RO extensions are two classic methods to deal with uncertain parameters. SP is suitable to the cases when the distribution of the underlying uncertainty is fully known. Some classical SP studies can be found in Cui et al. (2010), Noyan et al. (2016), and Paul and Zhang (2019) for more details. While SP has been shown to be an effective approach, the optimizers' curse will occur if data sample was calibrated from a biased distribution. In contrast, RO does not rely on exact probability distributions. Instead, it uses uncertain sets to represent random parameters and focuses on the worst-case scenario. It is particularly favored in contexts where uncertainties are difficult to estimate or data is scarce and has received more attention recently. Baron et al. (2011) built a multi-period capacitated fixed charge (single-stage) robust location model. Gülpınar et al. (2013) proposed a tractable (single-stage) RO method to approximately solve stochastic facility location problem.

Whereas the single-stage RO method determines only here-and-now decisions which require the optimal solution is feasible for all possible choices simultaneously, the two-stage RO method allows wait-and-see decisions that can adapt to the realized observations, while it is often possible in practice to adjust the solution after the uncertainty has been revealed (Ben-Tal et al., 2004). Several two-stage RO optimization models have been proposed to design supply chain networks in the presence of disruptions. For instance, Cheng et al. (2021) considered a fixed-charge location problem that permits customer reassignment as a recourse decision, where both demand and the number of disrupted facilities are uncertain and represented by budgeted uncertainty sets. Rahmati et al. (2022) proposed a two-stage robust optimization or the uncapacitated hub location problem under demand uncertainty, where an uncertainty budget was used to control the level of conservatism. The study also showed that an accelerated Benders decomposition with Pareto-optimal cuts improved solution efficiency, while a larger uncertainty budget led to higher total cost and more hub establishments. As demonstrated by Ben-Tal et al. (2004), even a basic two-stage RO problem can be NP-hard. Some researchers proposed different solution algorithms including approximation algorithms (Feng et al., 2023), Benders cutting plane algorithm (Yin et al., 2024), and the C&CG method. Note that C&CG algorithm generally requires significantly fewer iterations than Benders decomposition algorithm to find the optimal solution and has efficiency in solving two-stage RO models.

As mentioned above, there have been numerous studies related to single level optimization including SP and RO methods, which illustrated the application scenarios of different methods. Focusing on the disruptions risk, two-stage RO has been considered extensively as it is capable of generating less conservative solutions (Chassein et al., 2018). However, traditional single-level optimization models typically assume that the system designer and the network user share aligned objectives or belong to the same decision-making entity. As a result, decisions are made centrally from a system-wide perspective, without accounting for the decentralized behaviors or diverse priorities of operational agents. This simplification may overlook critical conflicts and behavioral heterogeneity in real-world supply chain systems, potentially leading to unintended spillover effects on the economy and environment. Given the long-term and wide-ranging impact of network design and facility location decisions, it is essential to incorporate the decentralized nature of decision-making and the potential conflicts between strategic planners and operational users.

2.2 Bilevel facility location

BO, which reflects the Stackelberg leader-follower game, dates back to the early publications by Bracken and McGill (1973). This methodology is usually used to formulate a systematic optimization problem with a typical leader-follower hierarchical structure. Recent developments in BO under uncertainty have been reviewed by Beck et al. (2023). A series of network design and facility location have been studied using BO schemes. For example, a two-stage stochastic BO model was applied to the established timberlands supply chain for biorefinery investments, where harvester and manufacturer are the individual decision makers (Yeh et al., 2015). Li et al. (2021) studied a facility location problem considering the boundedly rational attacker and probabilistic fortification.

Considering the facility disruption is a kind of risk with high impact but low frequency, a specific class of bilevel models, known as facility fortification models, has been developed to address such problems (Scaparra and Church, 2008). These models follow worst-case analysis and are typically formulated within a defend-attack framework, enabling the system designer to proactively safeguard against the most critical disruption scenarios, even without sufficient historical data. The defender, as the leader, selects facilities to protect in order to minimize the worst-case loss of efficiency, while the attacker, as the follower, tries to incur worst-case disruption impact by targeting unprotected facilities. The facility fortification problem, as a special class of BO models, has been extensively studied under uncertainty (Zhang & Özaltın, 2021; Zhang et al., 2022). We mention that other types of BO FLP models are much less studied.

However, compared to SP approaches, RO methods for handling uncertainty in BO models have rarely investigated (Beck et al., 2023). Traditional RO formulations, such as the classic 'min–max' structure, adds another level formulation on top of the underlying BO models. This significantly increases the computational complexity and often makes the robust BO models intractable. Xiong et al. (2021) proposed a robust bilevel resource recovery planning problem in the setting of the public and private partnership considering the uncertain feedstock condition, and reformulate equivalently to a mixed-integer linear program of reasonable size. Li et al. (2023) made a robust bilevel network design and pricing model considering yield uncertainty in biomass supply chain network, and converted it to a single level optimization. Overall, the existing robust bilevel studies are based on particular model reformulations, then convert to a solvable form, which may require special structures of the original models. To the best of our knowledge, apart from Xu (2023), no general solution algorithms exist that can exactly solve two-stage robust bilevel optimization models. Xu (2023) proposed C&CG-BO that explicitly handles general two-stage robust bilevel optimization. Nevertheless, we note that the standard C&CG-BO's is challenging for computational time. To do so, leveraging FLP's structures and properties to address practical-scale instances can effectively enhance the computation. Building on this foundation, we develop an improved algorithm of C&CG-BO, termed C&CG-BO-DDU, that exploit the structures and properties of 2-stage RBO-FLP to effectively address practical-scale instances.

**Table 1** Summary of related literature on facility location problem in facility location

| Authors | Problem | Uncertainty | Level | Stage | SP/RO | Solution method |
|---|---|---|---|---|---|---|
| Snyder and Daskin (2005) | P-median | Supply | Single | One | SP | Lagrangian relaxation |
| Paul and Zhang (2019) | Fixed-cost location | Supply | Single | Two | SP | MIP |
| Gülpınar et al., (2013) | Fixed-cost location | Demand | Single | One | RO | MIP |
| An et al. (2014) | P-median | Supply | Single | Two | RO | C&CG |
| Cheng et al. (2021a) | Fixed-cost location | Supply | Single | Two | RO | C&CG, affine policy |
| Hansen et al. (2004) | Fixed-cost location | / | Bilevel | / | / | Reformulation |
| Lamontagne et al. (2023) | Fixed-cost location | / | Bilevel | / | / | MIP |
| Yeh et al., (2015) | Supply chain allocation | Supply | Bilevel | Two | SP | Scenario programming |
| Li et al. (2021) | Fixed-cost location | Supply | Bilevel | Two | SP | Hybrid genetic algorithm |
| Xiong, et al., (2021) | Resource recovery plan | Supply | Bilevel | Two | RO | Reformulation |
| Li et al., (2023) | Network design and pricing | Supply | Bilevel | One | RO | KKT |
| This paper | Fixed-cost location | Supply | Bilevel | Two | RO | C&CG framework |

We summarized the features considered in the literature on the facility location problem in Table 1, from the perspectives of SP, RO, BO and solution algorithm. In summary, existing studies primarily have the following issues: 1) RO and SP are widely used for designing reliable and resilient supply chain networks, but they typically assume a single decision-maker by merging the roles of network designer and operator. This simplification overlooks the decentralized nature of real-world supply chains, where strategic planners and operational users may have distinct objectives and concerns. 2) BO provides a natural way to capture the hierarchical decision-making process, but its integration with robust methods to address low-probability, high-impact disruptions remains largely unexplored. In particular, two-stage robust bilevel models that incorporate facility disruption risk and decentralized operational behavior have received little attention.

To address the above issues, this paper proposes an integrated two-stage RO and BO framework that models the strategic facility location decisions of the designer and the responsive allocation behaviors of the operator, under disruption risk scenarios. Our work distinguishes itself by incorporating a decentralized decision-making perspective, explicitly accounting for the behavior of the network user within the context of robust facility location. To the best of our knowledge, this paper presents both the first model and exact solution algorithm(s) for two-stage robust bilevel facility location problems.

## 3. Model Formulation

3.1 Bilevel facility location formulation under facility disruption

Consider a two-echelon supply chain network composed of multiple capacitated facilities that serve a set of customers with a single product. Let $N$ denote the set of customers, and assume that the demand for each customer $i \in N$ is known by decision makers. The customers' demands are met by operable facilities, with each facility $j \in F$ having a maximum capacity of $K_j$. Let $c_{ij}$ represent the unit assignment cost to customer $i \in N$ from supply facility $j \in F$. Any unmet demand due to insufficient supply capacity is treated as a loss in the system and incurs a unit penalty $\rho_i$ for lost sales. We introduce $u_i$ auxiliary variable to represent the unsatisfied demand at customer node $i \in N$.

We consider two decision makers: the network designer and the network user with their distinguished pursuit in the system. The network designer first makes the strategic location decision, aiming to optimize long-term performance metrics such as cost efficiency. Once the facilities are established, the network user responds by making supply allocation decisions with the goal of minimizing the total unmet demand, thereby ensuring a high service level. Thus, the whole decision-making process between the network designer (leader) and the network user (follower) represents a Stackelberg game, and can be formulated as a BO model. In the upper-level problem, the network designer determines the optimal location to minimize the system total cost, including location cost, product allocation cost, and unmet demand penalty. In the lower-level, the network user determines the optimal allocation decisions to minimize the total unmet demand.

Recently, the increasing occurrence of black swan events has led to more frequent disruptions in supply facilities, resulting in significant uncertainty in facilities' supply capacities. As facility disruption is a low-probability, high-impact risk, and existing data is often insufficient to have a convincing probability distribution, we adopt two-stage RO approach, instead of the traditional stochastic optimization approaches, to better handle such uncertainty. The complete decision-making sequence and the involvement of disruptions is shown in Fig.2 below:

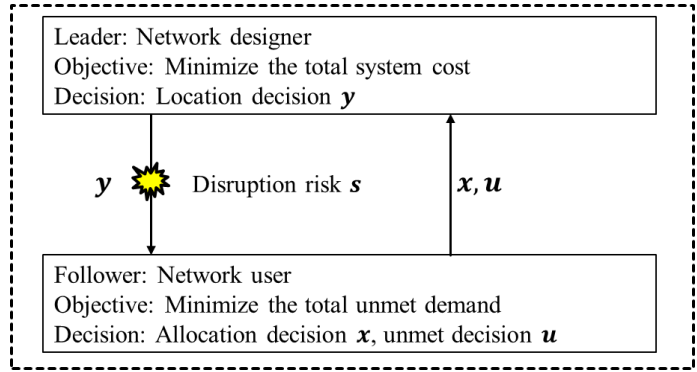

Fig.2 The framework of two-stage robust bilevel facility location under facility disruption

As shown in Fig.2, the network designer (leader) is a here-and-now decision maker, making location decisions before the risk occurs. Once disruptions occur, the network user (follower) responds as a wait-and-see decision maker, optimizing allocation and unmet demand decisions. Unlike the existing two-stage RO models, the network user at the operational level is different from the strategic designer, and they have distinguished benefit concerns. Hence, when minimizing the total system cost, the network designer now must account not only for facility disruption risks but also for the network user's behavior. Next, we present the concrete mathematical model of our two-stage robust bilevel model for FLP, which is referred to as 2-stage RBO-FLP hereafter in this paper.

Parameters

$f_j = $ fixed cost to open a facility $j \in F$

$c_{ij} = $ unit allocation cost on arc $(i,j) \in \mathcal{A}$

$K_j = $ maximum supply capacity of facility $j \in F$

$\rho_i = $ unit penalty cost for unmet demand at customer node $i \in N$

$d_i = $ demand quantity of customer node $i \in N$

$s_j = 1$ if facility $j \in F$ is disrupted, $s_j = 0$ otherwise

Decision variables

$y_j = 1$ if facility $j \in F$ is opened, $y_j = 0$ otherwise

$x_{ij} = $ product allocation decision on arc $(i,j)$

Auxiliary variable

$u_i = $ unmet demand at node $i \in N$

2-stage RBO-FLP

$$W_{\{RBO\}} = \min_{y \in \{0,1\}^{|F|}} \sum_{j \in F} f_j y_j + \max_{s \in \mathcal{U}} \min_{(x,u) \in \phi(y,s)} \left( \sum_{i \in N} \sum_{j \in F} c_{ij} x_{ij} + \sum_{i \in N} \rho_i u_i \right) \quad (1a)$$

Where $\phi(y,s)$ is the optimal solution set of the following lower-level problem

$$\phi(y,s) = \min_{x,u} \sum_{i \in N} u_i \quad (1b)$$

$$\sum_{i \in N} x_{ij} \leq K_j y_j (1 - s_j) \quad \forall j \in F \quad (1c)$$

$$\sum_{j \in F} x_{ij} + u_i = d_i \quad \forall i \in N \quad (1d)$$

$$x_{ij} \geq 0 \quad \forall i \in N, j \in F \quad (1e)$$

$$u_i \geq 0 \quad \forall i \in N \quad (1f)$$

$$\mathcal{U} = \left\{ s \,\middle|\, \begin{array}{l} \sum_{j \in F} s_j \leq \Gamma \\ s_j \in \{0,1\} \quad \forall j \in F \end{array} \right\} \quad (1g)$$

The upper-level objective (1a) is to minimize the system total cost including fixed location cost and the worst-case disruption risk scenario allocation cost and unmet penalty cost considering lower-level's behavior. The lower-level objective (1b) aims to minimize the total unmet demand to guarantee the service level after disruption risk scenarios are realized. Constraints (1c) model the facility capacity constraint after the disruption risk. Constraints (1d) mean that all demand is considered, whether allocated to a facility or unmet. Constraint (1g) in the budget uncertain set models the disruption risk level, denoting the total number of facilities being disrupted. The other constraints present individual variables restrictions. To facilitate our following discussion, we let $X(y,s)$ denote the feasible set of the lower-level problem, i.e., the set defined by (1b) - (1f).

As can be seen in Table 1, FPL has been studied using two-stage RO scheme, where first-stage location decisions are made prior to the realization of disruption and second-stage allocation decisions are determined after disruptions. As it is tightly related to 2-stage RBO-FLP, we next present the classical single level two-stage RO formulation, which is referred to as 2-stage RO-FLP.

$$W_{\{RO\}} = \min_{y \in \{0,1\}^{|F|}} \sum_{j \in F} f_j y_j + \max_{s \in \mathcal{U}} \min_{(x,u) \in X(y,s)} \left( \sum_{i \in N} \sum_{j \in F} c_{ij} x_{ij} + \sum_{i \in N} \rho_i u_i \right) \quad (1h)$$

3.2 Structural properties of 2-stage RBO-FLP

In the context of BO, one key concept is the high point problem, which simply drops off the lower-level objective function and merges the lower-level constraints with the upper-level ones to build an augmented single level formulation. Actually, this concept can be easily extended to the 2-stage RBO-FLP, and we note that the 2-stage RO-FLP model actually serves as the high-point problem of the former.

**Proposition 1.** (a) Any solution feasible to 2-stage RO-FLP is feasible to 2-stage RBO-FLP too. (b) 2-stage RO-FLP a relaxation to 2-stage RBO-FLP, and yields a lower bound for the 2-stage RBO-FLP model, i.e., $W_{\{RO\}} \leq W_{\{RBO\}}$.

**Proof.** Note that the second-stage lower-level problem in 2-stage RBO-FLP is feasible for any $(y,s)$. Also, 2-stage RBO-FLP and 2-stage RO-FLP share the same feasible set for the first stage and the uncertainty set. Hence, the first statement is naturally valid. Regarding the second statement, given that $\phi(y,s) \subseteq X(y,s)$ for arbitrary $(y,s)$, it simply follows. □

**Remark 1.** With Proposition 1, it is clear that, under the two-stage RO context, 2-stage RO-FLP provides a simpler basis to appreciate 2-stage RBO-FLP, which is analogous to what we have seen in the literature of deterministic BO. Hence, we note that it serves the high-point problem of 2-stage RBO-FLP. There is also a distinct feature of 2-stage RBO-FLP worth highlighting. Since the second-stage decision is completely determined by the follower in 2-stage RBO-FLP, instead of saying that the second-stage has the **relatively complete recourse property**, we say the follower has the **relatively complete responsiveness property**, given that the lower-level problem has a finite optimal value for any $y, s$ feasible to the first stage and belonging to $\mathcal{U}$.

**Remark 2.** If $c_{ij} = 0, \rho_i = 1, \forall i \in N, j \in F$, 2-stage RBO-FLP model is also a special two-stage robust single level optimization model. In this case, the lower-level and upper-level optimal solutions coincide. This observation can be generalized to the following result.

**Proposition 3.** If $\rho_i = \rho \geq \max\{c_{ij}, \forall i \in N, j \in F\}$, we have $W_{\{RBO\}} = W_{\{RO\}}$.

**Proof.** Consider a fixed 2-tuple $(y^*, s^*)$. There are two cases according to the total capacity from survived facilities.

In the first case, we have enough total capacity to serve customers. It follows that for both 2-stage RBO-FLP and 2-stage RO-FLP, we have $\sum_{i \in N} u_i^* = 0$, which renders $u_i^* = 0$ for all customers in any optimal solution to the recourse or the lower-level problem, respectively. Suppose that $(x^*, \mathbf{0}) \in X(y^*, s^*)$ is an optimal solution to the recourse problem of 2-stage RO-FLP. It is clearly optimal to the lower-level problem of 2-stage RBO-FLP, i.e., $(x^*, u^*) \in \phi(y^*, s^*)$. Let $W_{\{RO\}}(y^*, s^*)$ and $W_{\{RBO\}}(y^*, s^*)$ denote the second stage costs in 2-stage RBO-FLP and 2-stage RO-FLP respectively. We have $W_{\{RO\}}(y^*, s^*) \leq W_{\{RBO\}}(y^*, s^*) \leq W_{\{RO\}}(y^*, s^*)$.

In the second case, the total capacity from survived facilities is not sufficient to serve all customers. Let $\widehat{F}(y^*, s^*)$ denote the set of survived facilities. Suppose that $\left(x^*, \sum_{\{i \in N\}} d_i - \sum_{\{j \in \widehat{F}(y^*, s^*)\}} K_j\right) \in X(y^*, s^*)$ is an optimal solution to the recourse problem of 2-stage RO-FLP.

Again, it is optimal to the lower-level problem of 2-stage RBO-FLP, i.e., $\left(x^*, \sum_{\{i \in N\}} d_i - \sum_{\{j \in \hat{F}(y^*, s^*)\}} K_j\right) \in X(y^*, s^*)$. By using the same argument to the first case, we have $W_{\{RO\}}(y^*, s^*) \leq W_{\{RBO\}}(y^*, s^*) \leq W_{\{RO\}}(y^*, s^*)$.

For any $(y^*, s^*)$, the conclusion follows directly. □

This implies that a high unmet penalty encourages or facilitates full cooperation between the network designer and the network user. Next, we illustrate the differences in $(x, u)$ between 2-stage RO-FLP and 2-stage RBO-FLP under a given location decision and disruption scenario $(y, s)$ through an example to demonstrate that the optimal solution obtained from the 2-stage RBO-FLP model fully utilizes the available capacity. As depicted in Fig. 3, the diamond represents a supply facility, with the number inside indicating its maximum available supply capacity. The green area within the diamond represents the actual used supply capacity. The circle represents a demand point, where the number inside denotes the demand quantity, and the yellow area within the circle represents the unmet demand. The transportation cost from the supply facility to the three demand points are 1, 1, and 1.41, respectively. The red lines represent allocation decisions in the 2-stage RO-FLP model, while the blue lines represent allocation decisions in the 2-stage RBO-FLP model.

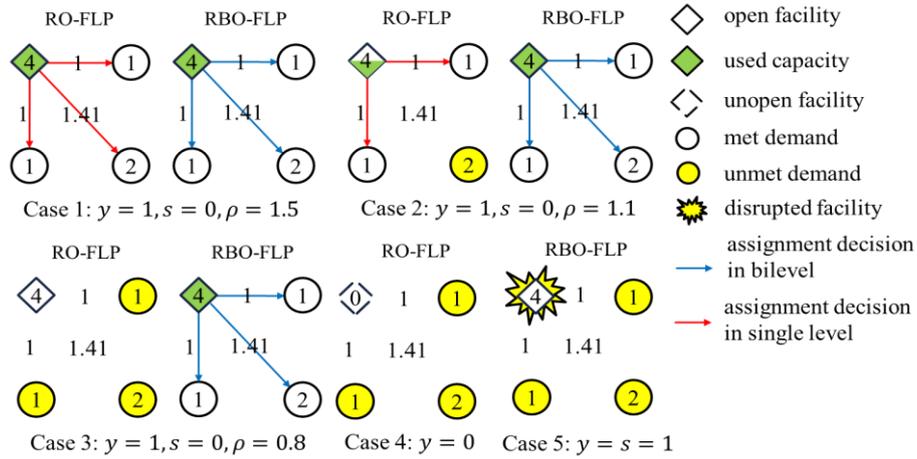

Fig.3  The differences in $(x, u)$ between 2-stage RO-FLP and 2-stage RBO-FLP

Cases 1, 2, and 3 illustrate the impact of unmet penalty on $(x, u)$. In the 2-stage RO-FLP model, if the transportation cost exceeds the unmet penalty, the corresponding demand will not be fulfilled, leading to wasted supply capacity (as indicated by the green area in the diamond being smaller than the total diamond area). However, in the 2-stage RBO-FLP model, since the network users aim to minimize total unmet demand, the network designer fully uses the available supply capacity to meet demand (as indicated by the green area in the diamond being equal to the total diamond area). The 2-stage RBO-FLP model better reflects real-world scenarios in which ensuring supply after the disruption risk is critical, highlighting the necessity of a bilevel modeling approach. Cases 4 and 5 demonstrate that when no facilities are operational, the allocation decisions $(x, u)$ remain identical in both models.

In the following, we consider another simple case where no facility will be built in an optimal solution. It can be easily proven by arguing that letting all demand unmet is actually better than building and utilizing any facility.

**Proposition 4.** If $\rho_i = \rho \leq \min \{c_{ij}, \forall i \in N, j \in F\}$, i.e., the penalty for unmet demand is identical across all demand nodes and is lower than the minimum unit transportation cost. The optimal solutions $(y, x, u)$ of the 2-stage RBO-FLP model and the 2-stage RO-FLP are both $(\mathbf{0}, \mathbf{0}, \mathbf{d})$.

Combining Propositions 3 and 4, it can be seen that a low unmet demand penalty discourages the cooperation between the network designer and the network user, and a high penalty promotes such cooperation. In the next section, we develop an improved and more efficient algorithm, C&CG-BO-DDU, to solve 2-stage FLP-RBO that drastically improves its computational performance.

## 4. Solution methodology based on C&CG framework

Compared with the standard 'min-max-min' structure of 2-stage RO-FLP model, 2-stage RBO-FLP model is more challenging to solve due to the additional lower-level 'min' problem embedded in the second stage. Lots of literature has proven the effectiveness of (basic) C&CG algorithm in solving two-stage robust RO problems (Zeng and Zhao, 2013; Cheng et al., 2021; Wang and Qi, 2020). Xu (2023) employed the C&CG framework and develop C&CG-BO to solve general two-stage robust bilevel problems. Similar to basic C&CG solving the 2-stage RO-FLP model (See Appendix A), C&CG-BO involves solving the master problem (MP) and the subproblem (SP). Note that MP, representing a BO model built with respect to a set of scenarios, derives the optimal first-stage facility location decision, and the SP is solved to identify the worst facility disruption scenario. We highlight that, while the MP is a regular BO problem, SP is actually a pessimistic BO problem that is known much more difficult to solve (Wiesemann et.al, 2013; Zeng 2020). Actually, in our computationally study, solving SP is the bottleneck of the whole solution process. Nevertheless, by leveraging a critical structure underlying the DDU set, we develop C&CG-BO-DDU, which allows the SP to be solved drastically faster and, in turn, improve our solution capacity on large-scale 2-stage FLP-RBO instances by an order of magnitude.

4.1 MP of C&CG-BO

For a set of given disruption scenarios $s^{l*}$, MP of C&CG-BO is a BO model as in the following:

MP

$$\min \eta \tag{2a}$$

$$\eta \geq \sum_{j \in F} f_j y_j + \sum_{i \in N} \sum_{j \in F} c_{ij} x_{ij}^l + \sum_{i \in N} \rho_i u_i^l \quad l = 1,2,3 \ldots, \mathcal{L} \tag{2b}$$

$$(x^l, u^l) \in \phi^l(y, s^{l*}) \quad \text{with} \tag{2c}$$

$$\phi^l(y, s^{l*}) = \text{argmin} \left\{ \begin{array}{l} \sum_{i \in N} u_i^l : \sum_{i \in N} x_{ij}^l \leq K_j y_j (1 - s_j^{l*}) \quad \forall j \in F; \\ \sum_{j \in F} x_{ij}^l + u_i^l = d_i \quad \forall i \in N; \\ x^l \geq 0, u^l \geq 0 \end{array} \right\}; l = 1,2,3 \ldots, \mathcal{L}$$

Note that allocation variables $x_{ij}^l$, auxiliary variables $u_i^l$ and $\phi^l(y, s^{l*})$ are defined with respect to $l$, indicating their associations with $s_j^{l*}$, the $l$-th worst-case realization of uncertain variables $s_j$ obtained by solving the corresponding SP in the $l$-th iteration. To convert it into a format that is readily computable for a solver, $\phi^l(y, s^{l*})$ can be replaced by the corresponding optimality conditions, e.g., KKT conditions, and then reformulated by linearization techniques to derive an equivalent MIP form of MP. This equivalent MIP formulation can be seen in Appendix B.1. Since the given set of disruption scenarios is a subset of $\mathcal{U}$, MP is naturally a relaxation of 2-stage RBO-FLP and provides a lower bound.

4.2 SP of C&CG-BO

SP is solved to identify the worst disruption risk scenarios under given $\mathbf{y}$, which, in its original form, is written as:

SP-1

$$\max_{s\in\mathcal{U}} \min_{x,u} \sum_j f_j y_j^* + \sum_i \sum_j c_{ij} x_{ij} + \sum_i \rho_i u_i \tag{3a}$$

$$(\mathbf{x}, \mathbf{u}) \in \phi(\mathbf{y}^*, \mathbf{s}) = \operatorname{argmin} \left\{ \begin{array}{l} \sum_i u_i : \sum_i x_{ij} \leq K_j y_j^*(1-s_j) \quad \forall j \in F \\ \sum_j x_{ij} + u_i = d_i \quad \forall i \in N \\ \mathbf{x} \geq 0, \mathbf{u} \geq 0 \end{array} \right\} \tag{3b}$$

As mentioned, SP-1 is a challenging pessimistic BO model that actually has a three-level structure. Note that its structure prevents us from directly applying any optimality conditions to transform it into a single level MIP. Zeng (2020) proposed a level reduction technique to address such a challenging structure, which allows us to develop a simpler MIP equivalence of SP-1 that is readily computable.

**Proposition 5**. The pessimistic BO model SP-1 can be equivalently transferred to the following BO model (referred to as SP-BO).

SP-BO

$$\max_{s\in\mathcal{U}} \sum_j f_j y_j^* + \sum_i \sum_j c_{ij} \overline{x}_{ij} + \sum_i \rho_i \overline{u}_i \tag{4a}$$

$$\sum_i \overline{x}_{ij} \leq K_j y_j^*(1-s_j) \quad \forall j \in F \tag{4b}$$

$$\sum_j \overline{x}_{ij} + \overline{u}_i = d_i \quad \forall i \in N \tag{4c}$$

$$\overline{x} \geq 0, \overline{u} \geq 0 \tag{4d}$$

$$(\mathbf{x}, \mathbf{u}) \in \tilde{\phi}(\mathbf{y}^*, \mathbf{s}, \overline{\mathbf{x}}, \overline{\mathbf{u}}) = \operatorname{argmin} \left\{ \begin{array}{l} \sum_j f_j y_j^* + \sum_i \sum_j c_{ij} x_{ij} + \sum_i \rho_i u_i : \\ \sum_i x_{ij} \leq K_j y_j^*(1-s_j) \quad \forall j \in F \\ \sum_j x_{ij} + u_i = d_i \quad \forall i \in N \\ \sum_i u_i \leq \sum_i \overline{u}_i \\ \mathbf{x} \geq 0, \mathbf{u} \geq 0 \end{array} \right\} \tag{4e}$$

**Proof.** See Appendix C.

Again, SP-BO is a regular BO model that can be reformulated as an MIP, referred to as SP-MIP, which is presented in Appendix B.2. Note that solving SP-MIP identifies the worst-case facility disruption scenario based on the given first-stage location decisions, thereby providing a valid upper bound for the 2-stage RBO-FLP model.

Noting that $\mathcal{U}$ is a finite set and every iteration SP-MIP generates a new scenario, the next result regarding the convergence of C&CG-BO follows directly.

**Proposition 6**. Algorithm 1 generates an optimal solution within a finite number of iterations.

4.3 Algorithm procedures with an enhancement strategy

Although C&CG-BO can converge within a small number of iterations, we observe that solving SP-MIP, the reformulation of the pessimistic bilevel subproblem, could be extremely time-consuming, which leads to a heavy computational burden in solving 2-stage RBO-FLP. Fortunately, by exploiting the DDU structure of $\mathcal{U}$, C&CG-BO-DDU has a significant and critical impact on our solution approach. Recall that the assumption made in Section 3 can be directly captured by the following constraints.

$$s_j \leq y_j \quad \forall j \in F \tag{4f}$$

By incorporating (4f) into $\mathcal{U}$ and revising the lower-level problem, we reformulate and have an alternative form of SP-1.

SP-2

$$\max_{s \in \mathcal{U}(y)} \min_{x,u} \sum_j f_j y_j^* + \sum_i \sum_j c_{ij} x_{ij} + \sum_i \rho_i u_i$$

$$(x, u) \in \phi(y^*, s) = \operatorname{argmin} \left\{ \sum_i u_i : \begin{array}{l} \sum_i x_{ij} \leq K_j y_j^*(1 - s_j) \quad \forall j \in F \\ \sum_i x_{ij} \leq K_j y_j^* \quad \forall j \in F \\ \sum_j x_{ij} + u_i = d_i \quad \forall i \in N \\ x \geq 0, u \geq 0 \end{array} \right\}$$

$$\mathcal{U}(y) = \left\{ s \;\middle|\; \begin{array}{l} \sum_{j \in F} s_j \leq \Gamma \\ s_j \leq y_j \quad \forall j \in F \\ s_j \in \{0,1\} \quad \forall j \in F \end{array} \right\}$$

**Proposition 7.** (a) $\mathcal{U}(y)$ is a subset of $\mathcal{U}$. (b) SP-1 and SP-2 have the same optimal value.

**Proof.** As the first statement is clear, we consider the second statement only. We first show that the optimal value of SP-1 is greater than or equal to that of SP-2. This is because SP-2 contains additional constraints of the form $\sum_i x_{ij} \leq K_j y_j^*$ $\forall j \in F$, which lead to a smaller feasible region and SP-2 is solved over a subset of the feasible set of SP-1.

Next, we show that the optimal value of SP-2 is greater than or equal to that of SP-1. Consider the optimal solution $(x^*, u^*)$ of SP-1 also satisfies the additional constraints in SP-2 — namely $\sum_i x_{ij} \leq K_j y_j^*$ $\forall j \in F$ and $s_j \leq y_j$ $\forall j \in F$ — it is feasible in SP-2 as well. Combining both directions, we conclude that SP-1 and SP-2 have the same optimal value. □

**Remark 3:** (1) We highlight that, after incorporating (4f), $\mathcal{U}$ is converted into a DDU set, which, according to Proposition 7, helps us reduce the search space of $s$. In our numerical experiments, such reduction has a great and positive impact on solving SP. Actually, the computation can speed up by more than one order of magnitude, which drastically improves our solution capacity to the large and challenging SP problems. Regarding MP in section 4.1, we do not make any changes. (2) Such a great computational performance enhancement verifies and strengthens one observation made in Zeng and Wang (2022) that converting an uncertain set into a DDU set could advance our computational capacity significantly. Hence, we do believe that such conversion and related techniques are important future research topics.

4.4 Algorithm Framework

In summary of the above proposition and analysis, our C&CG-BO-DDU for two-stage RBO-FLP is shown in Fig. 4 and the detailed procedure is given in Algorithm 1.

| Algorithm 1: C&CG-BO-DDU algorithm for the 2-stage RBO-FLP model | |
|---|---|
| Step1 | Initialization: Let lower bound $LB = -\infty$, upper bound $UB = +\infty$, iterations $l = 0$, $s^{0*} = 0$, $\epsilon = 0.1\%$. |
| Step2 | Solve the MP-MIP to obtain optimal location decision $y^*$. Set $LB = max\{LB, \eta\}$. |
| Step3 | Solve the SP-MIP to find the worst-case disruption risk scenario $s^{l*}$. Let $\psi$ be the objective value of the SP. Set $UB = min\{\psi, UB\}$. |
| Step4 | If $Gap = \frac{UB-LB}{UB} \leq \epsilon$, the algorithm terminates; else, set $l = l + 1$, add the identified worst-case disruption risk scenario $s^{l*}$ and its corresponding recourse decision variables and constraints to the MP-MIP and go to Step2. |

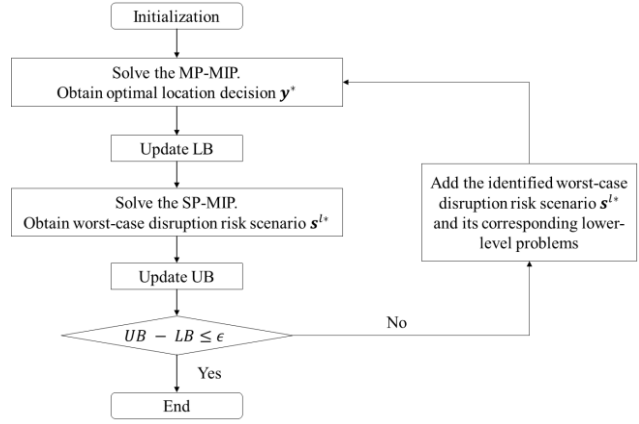

Fig.4 The C&CG-BO-DDU algorithm for the 2-stage RBO-FLP model

## 5. Numerical study and analysis

In this section, we conduct computational experiments to evaluate the performance of the proposed 2-stage RBO-FLP model in comparison with the traditional 2-stage RO-FLP model. Specifically, we compare key metrics including total cost, service level, unit service cost, and average supply capacity utilization. Finally, we further analyze the conditions under which centralized or decentralized decision-making structures are more appropriate. The algorithm is implemented in Python and executed on a PC equipped with an Intel i7 2.60 GHz processor, 16 GB RAM, and a 64-bit Windows 10 operating system. All models are solved using Gurobi 9.0.

5.1 Generation of random data sets and basic experimental setting

The test dataset is constructed based on a primary dataset of 49 points (Daskin, 2011). To enhance diversity, we generate an additional dataset of the same size and merge them. A subset of 46 points is then randomly selected for testing, ensuring that the attributes of the generated data remain within the range of the original dataset. From the test dataset, six points are randomly designated as potential supply facilities, while the remaining points serve as customer nodes. The other model parameters are detailed in Table 2.

Table 2 Template for the values of parameters

| Parameters | Values | Parameters | Values |
|---|---|---|---|
| $d_i$ | the state population $*10^{-4}$ | $\rho_i$ | $0.01 * \sum_{j \in F} f_j / |F|$ |

| $f_j$ | median home value in the city $*10^{-2}$ | $K_j$ | $1.2 * \sum_{i \in N} d_i / |F|$ |
|---|---|---|---|
| $c_{ij}$ | the distance between nodes $i$ and $j$ | | |

5.2 The performance of the C&CG algorithm

In this section, we show the performance of the proposed C&CG algorithm, and specifically compare the basic C&CG-BO algorithm with its improved version, C&CG-BO-DDU, under different disruption risk levels $\Gamma$. To do so, we use $|F|$ and $|N|$ to represent the cardinality of sets $F$ and $N$, respectively. The column "Ratio" indicates the relative difference in computation time between the basic C&CG-BO and C&CG-BO-DDU algorithms, with smaller values meaning faster solving time for C&CG-BO-DDU. As shown in Table 3, C&CG-BO-DDU generally reduces and computation time compared to the basic C&CG-BO algorithm. As expected, solving time increases as $|F|$ grows, highlighting the need for scalable approximation algorithms for large-scale scenarios in future research.

Table 3 Computational Performance for algorithm

| | $(|F|, |N|) = (6,40)$ | | | $(|F|, |N|) = (10,40)$ | | | $(|F|, |N|) = (20,30)$ | | |
|---|---|---|---|---|---|---|---|---|---|
| $\Gamma$ | C&CG-BO Time (sec) | C&CG-BO-DDU Time (sec) | Ratio | C&CG-BO Time (sec) | C&CG-BO-DDU Time (sec) | Ratio | C&CG-BO Time (sec) | C&CG-BO-DDU Time (sec) | Ratio |
| 1 | 15 | 25 | 63% | 5921 | 3995 | -33% | 18050 | 19016 | 5% |
| 2 | 72 | 16 | -77% | 7750 | 4582 | -41% | 51909 | 20876 | -60% |
| 3 | 39 | 15 | -62% | 11058 | 4356 | -61% | 148410 | 163194 | 10% |
| 4 | 46 | 13 | -72% | 23859 | 6067 | -75% | 74279 | 60910 | -18% |
| 5 | 48 | 15 | -69% | 44749 | 2058 | -95% | 271711 | 21214 | -92% |

5.3 The impact of disruption risk level on the 2-stage RBO-FLP & 2-stage RO-FLP

In this section, we analyze the impact of the facility disruption risk level $\Gamma$ on the total system cost and service level of the 2-stage RBO-FLP in comparison with 2-stage RO-FLP models. Fig. 5 illustrates the total service volume and total cost of two models under increasing levels of disruption risk, denoted by markers 0 to 6. As observed in Fig. 5, the total cost for both models increase with a higher disruption risk level. However, the red line consistently lies below the blue line, indicating that the total system cost of the 2-stage RO-FLP model never exceeds that of the 2-stage RBO-FLP model. This result confirms that the 2-stage RO-FLP model serves as a relaxation of the 2-stage RBO-FLP model, providing a lower bound for its total cost. In real-world scenarios, the network designer of the 2-stage RBO-FLP model must account not only for minimizing total system cost but also for reducing overall unmet demand at the operational level, leading to additional cost considerations.

Although the 2-stage RBO-FLP model incurs a higher total cost, however, it achieves a higher level of service coverage. To further illustrate the advantage of the 2-stage RBO-FLP model, we introduce the concept of unit service cost (USC) defined as the cost for meeting unit service to reflects the service efficiency. A lower unit service cost indicates higher efficiency in resource utilization. We evaluate the unit service cost of both models under disruption risk levels 0 to 5 (risk level 6 is excluded as no facilities are open, resulting in a total assignment of zero). The results are shown in Fig. 6, which demonstrates that the 2-stage RBO-FLP model generally provides a more efficient level of service.

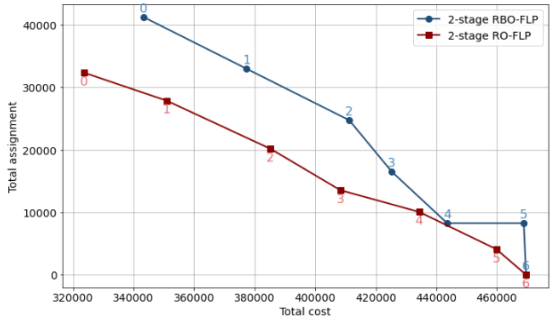 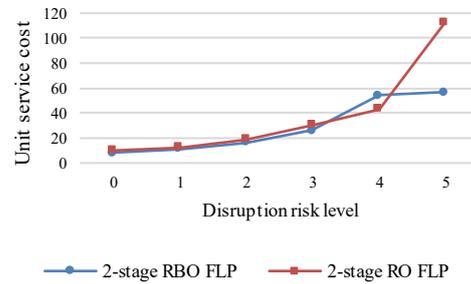

Fig.5　The impact of the disruption risk level　　　Fig.6　The unit service cost of the 2-stage RO-FLP model

To reflect the condition under which the two models apply, we introduce the concept of cost ratio and service ratio, of the model 2-stage RBO-FLP to model 2-stage RO-FLP. The experiment results are shown in Fig. 7, of which the x-axis represents the cost ratio, while the y-axis shows the service ratio. Each point corresponds to a risk level and the red dashed line indicates the threshold where the increase in cost is exactly matched by the increase in fulfilled demand. Points located above the red line indicate that the 2-stage RBO-FLP model achieves a greater improvement in service than that in cost. In such cases, the decentralized structure decision model, i.e. 2-stage RBO-FLP leads to higher service efficiency compared to the centralized structure decision model, 2-stage RO-FLP. Conversely, points below the line suggest the centralized strategy of the 2-stage RO-FLP model offers better cost-effectiveness. Most of the points lie above this reference line, indicating that the 2-stage RBO-FLP model often delivers a greater proportional improvement in service than total cost. This highlights the necessity of adopting the 2-stage RBO-FLP model in many real-world scenarios, where a moderate increase in investment can lead to significantly improved service guarantees.

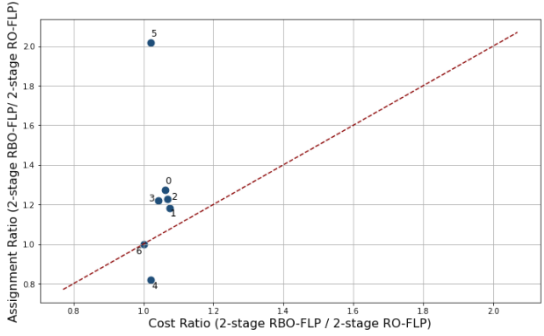
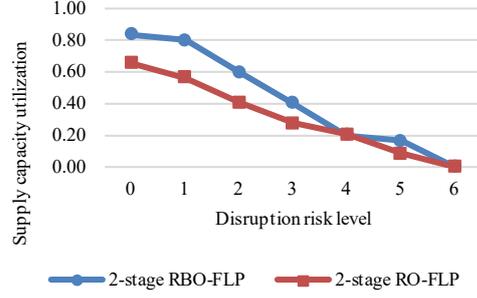

Fig.7  Cost and assignment ratio          Fig.8  Supply capacity utilization

5.4 The analysis of the average supply capacity utilization

In this section, we analyze the effective utilization of available supply capacity to validate the results of Fig. 3. Specifically, we define average supply capacity utilization as the ratio between the total actually allocated supply capacity and the total available supply capacity across all facilities. This metric reflects how efficiently the model makes use of the system's potential supply capacity. The average supply capacity utilization can be influenced by disruption risk (as seen in Case 5 of Fig. 3) or by insufficient allocation (corresponding to the 2-stage RO-FLP results in Cases 2 and 3 of Fig. 3). The detailed formulation of the average supply capacity utilization $\omega$ can be seen as follows:

$$\omega = \sum_{j \in F} \left( \frac{\sum_{i \in N} x_{ij}^*}{y_j^* * K_j} \right) / \sum_{j \in F} y_j^*$$

Where $y_j^*$ and $x_{ij}^*$ are the optimal solutions obtained by calculating the 2-stage RO-FLP model and the 2-stage RBO-FLP model respectively. We calculate the average supply capacity utilization $\omega$ of the above two models under different disruption levels, with the results shown in Fig. 8. The red curve, representing the results of the 2-stage RO-FLP model, is lower than the blue curve, which represents the results of the 2-stage RBO-FLP model. This indicates that the 2-stage RO-FLP model does not fully account for the concerns of the network user, leading to underutilization of available capacity. In contrast, the 2-stage RBO-FLP model exhibits higher average supply capacity utilization and thus provides a more realistic representation of emergency supply allocation after disruptions, better reflecting practical supply assurance needs.

Under a specific disruption risk level of $\Gamma = 2$, Fig. 9 presents the solutions of the 2-stage RO-FLP model and the 2-stage RBO-FLP model. The black dashed lines in Fig. 9 represent the specific allocation decisions between supply points (blue squares) and demand points (yellow circles). In Fig. 9a, the 2-stage RO-FLP model results in sparse allocations—primarily due to the high transportation cost relative to the unmet demand penalty, causing the network designer to sacrifice allocation coverage to minimize total cost. In contrast, Fig. 9b shows that the 2-stage RBO-FLP model achieves significantly higher allocation under the same number of available supply facilities. This highlights the advantage of the 2-stage RBO-FLP model by accounting for the service level of customers.

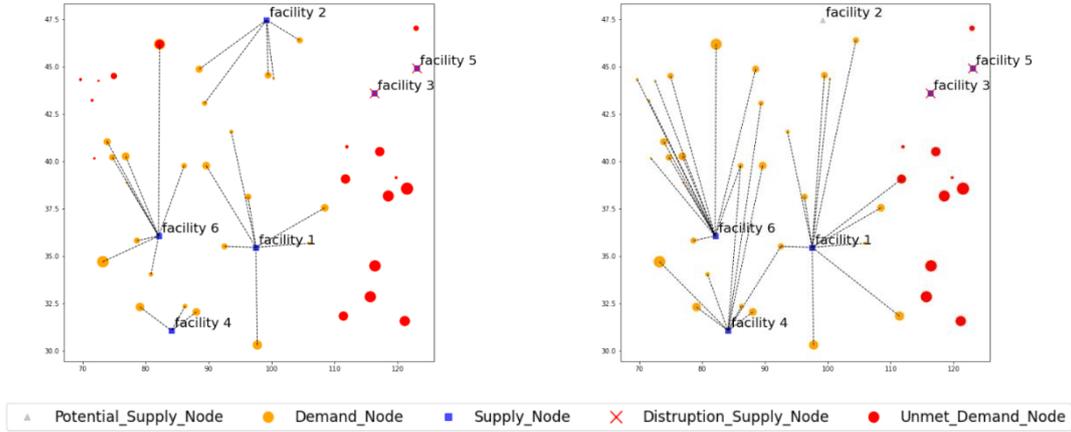

Fig.9a  The solution of the 2-stage RO-FLP model    Fig.9b  The solution of the 2-stage RBO-FLP model

Fig.9   The comparison of the solutions

5.5 The impact of the unmet demand penalty cost under varying disruption risk level

In Section 3, we present the Proposition 3 and Proposition 4 for analyzing the impact of the unmet demand penalty cost efficient $\rho_i$ on the 2-stage RO-FLP and 2-stage RBO-FLP models. In this section, we further investigate this effect under different facility disruption risk levels. We assume that the unmet demand penalty cost $\rho_i$ is uniform across all demand points, denoted as $\rho$, and set to the 0th, 25th, 50th, 75th, and 100th percentiles of $c_{ij}$ for all $i \in N, j \in F$, where the 0th and 100th percentiles correspond to the minimum and maximum values, respectively. For each disruption risk level and penalty cost scenario, we compute the total number of open facilities of the 2-stage RO-FLP and 2-stage RBO-FLP models, denoted as $\sum_j y_j^{S*}$ and $\sum_j y_j^{B*}$, as well as the total allocated capacity, denoted as $\sum_i \sum_j x_{ij}^{S*}$ and $\sum_i \sum_j x_{ij}^{B*}$. Here, the superscript $S$ represents results from the 2-stage RO-FLP model, while $B$ represents results from the 2-stage RBO-FLP model. The impact of the unmet demand penalty cost on 2-stage RBO-FLP and 2-stage RO-FLP models under varying disruption risk level are depicted in Fig. 10, where horizon axis and left vertical axis represent disruption risk level and penalty cost coefficient $\rho$ in percentile of transportation cost respectively.

Fig. 10a illustrates the difference $\sum_j y_j^{S*} - \sum_j y_j^{B*}$ under varying facility disruption risk levels $\Gamma$ and unmet demand penalty coefficient $\rho$ relative to transportation cost. Darker colors indicate larger differences, meaning the 2-stage RO-FLP model tends to open more facilities. White areas represent a difference of zero, implying that both models open the same number of facilities. As shown in Fig. 10a, the 2-stage RO-FLP model generally results in a higher number of open facilities, which aligns with the earlier average supply capacity utilization analysis and highlights a key distinction between the two models. In the 2-stage RBO-FLP model, the network designer accounts for user's decisions, leading to more conservative facility location decisions in response to disruption risks. Consequently, decentralized decision-making in the 2-stage RBO-FLP model tends to be more cautious compared to the centralized decision-making in the 2-stage RO-FLP model, resulting in fewer open facilities.

Fig. 10b illustrates the difference $\sum_i \sum_j x_{ij}^{S*} - \sum_i \sum_j x_{ij}^{B*}$ under varying facility disruption risk levels $\Gamma$ and unmet demand penalty efficient $\rho$. The red color spectrum indicates higher total allocation in the 2-stage RO-FLP model ($\sum_i \sum_j x_{ij}^{S*} \geq \sum_i \sum_j x_{ij}^{B*}$), while the blue spectrum represents higher total allocation in the 2-stage RBO-FLP model ($\sum_i \sum_j x_{ij}^{B*} \geq \sum_i \sum_j x_{ij}^{S*}$). Darker shades reflect larger differences between the two models. As shown in Fig. 10b, when the disruption risk level is relatively low ($\Gamma \leq 3$), the 2-stage RBO-FLP model exhibits higher total allocation. However, under high disruption risk levels and high unmet demand penalties (e.g., in the lower-right region of Fig. 10b), the 2-stage RO-FLP model results in greater total allocation. These findings highlight an important managerial insight: under low-risk conditions, upper-level network designers should fully account for lower-level network users' objectives in a decentralized way to ensure effective supply allocation in disruption scenarios. Conversely, when facing significant risks, lower-level network users should completely collaborate with upper-level network designers in a centralized way to mitigate risk more effectively. Finally, the first row of the results in Fig. 10a and Fig. 10b validates Proposition 4, confirming that when the unmet penalty is very low, the optimal solutions $(\boldsymbol{y}, \boldsymbol{x}, \boldsymbol{u})$ of the 2-stage RBO-FLP model and the 2-stage RO-FLP are both $(\boldsymbol{0}, \boldsymbol{0}, \boldsymbol{d})$. Consequently, $\sum_j y_j^{S*} - \sum_j y_j^{B*} = 0$ and $\sum_i \sum_j x_{ij}^{S*} - \sum_i \sum_j x_{ij}^{B*} = 0$.

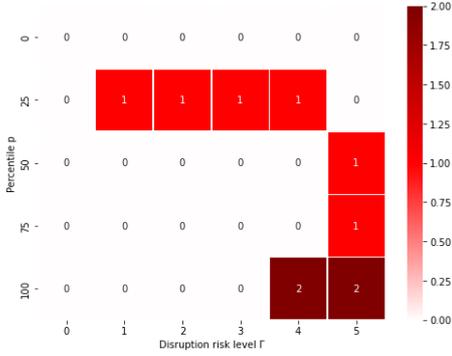
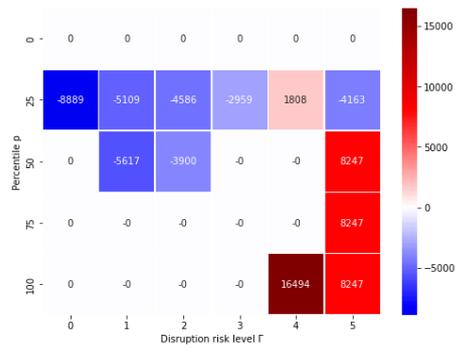

Fig.10a The analysis of location decisions    Fig.10b The analysis of allocation decisions

Fig.10 The impact of the unmet demand penalty cost under different risk level

## 6. Conclusion

In this paper, we present a 2-stage RBO-FLP model under facility disruption, by which, the network designer makes location decisions with minimum total cost before disruption risk, and network user makes optimal allocation and unmet decisions with minimum total unmet demand to guarantee service level after disruption risk occurs in a decentralized way. By utilizing C&CG-BO-DDU, we adopt an enhanced approach to address the complex SP by reducing the search space of the worst-case disruption scenario and improve the efficiency of solving the two-stage RBO-FLP drastically. The MP is a regular BO with multiple lower-level problems, while the SP is a pessimistic BO model. Through a series of reformulations, both the master and subproblems can ultimately be transformed into mixed-integer programming formulation that can be directly solved using existing optimization solvers. Through numerical experiments, we demonstrate that compared to the two-stage robust single-level optimization model (2-stage RO-FLP), the proposed 2-stage RBO-FLP achieves a better balance between cost investment and service level when facing facility disruption risk. Some theoretical properties are proposed and experiments are conducted to show the superiority of the proposed 2-stage RBO-FLP model in comparison with 2-stage RO-FLP model. The 2-stage RBO-FLP model extends both robust facility location problems and bilevel facility location problems and the enhancement strategy algorithm incorporating the DDU set can also be extended to improve the computational efficiency of other RO problems by reducing the search space of the uncertain scenarios.

Some managerial insights can be drawn as follows: 1) The baseline comparison between 2-stage RBO-FLP and 2-stage RO-FLP provides valuable guidance for managers in choosing between centralized and decentralized decision-making structures. By evaluating the trade-offs between total cost and service efficiency, managers can identify the more appropriate strategy under disruption risk. In many cases, the 2-stage RBO-FLP with decentralized decision-making leads to lower unit service costs, indicating higher service efficiency and better responsiveness to local demand conditions. 2) By incorporating the operational-level objective of minimizing total unmet demand, the network designer can improve average supply capacity utilization and better meet overall demand. 3) When facing lower levels of disruption risk and unmet demand penalty, adopting a decentralized decision-making structure can improve demand fulfillment by allowing the network designer to better accommodate the user's operational preferences. In contrast, under high disruption risk level and severe unmet demand penalties, a centralized approach is more appropriate, enabling closer coordination between the network designer and user to implement more effective risk mitigation strategies.